\documentclass[12pt,reqno]{amsart}
\usepackage{amssymb}
\usepackage{amsmath}

 \rm

\renewcommand{\(}{\left(}
\renewcommand{\)}{\right)}
\renewcommand{\[}{\left[}
\renewcommand{\]}{\right]}
\newcommand{\<}{\langle}
\renewcommand{\>}{\rangle}

\renewcommand{\bar}{\overline}
\newcommand{\abs}[1]{\left\lvert#1\right\rvert}
\newcommand{\norm}[1]{\left\lVert#1\right\rVert}
\newcommand{\st}{\:|\:}

\newcommand{\C}{{\mathbb{C}}}
\newcommand{\R}{{\mathbb{R}}}

\newcommand{\N}{{\mathbb{N}}}

\renewcommand{\phi}{\varphi}

\renewcommand{\H}{{\mathcal{H}}}

\theoremstyle{plain}
\newtheorem{thm}{Theorem}[section]

\theoremstyle{definition}

\theoremstyle{remark}
\newtheorem{rem}[thm]{Remark}

\title[pair correlation of Hermite coefficients]{On pair correlation of
  Hermite coefficients of functions from the Hardy class}

\author{Manish Chaurasia}
\address{Department of Mathematical Sciences, IIT (BHU), Varanasi 221005}
\subjclass[2010]{Primary 42C10; Secondary 42C05, 42B35, 33C45, 44A20}
\keywords{Bargmann transform, Hardy's theorem, Hermite function} 

\begin{document}

%% \begin{abstract}
%% Assuming that a function and its Fourier transform are dominated by 
%% centered Gaussians, in this article, we continue our study for the decay
%% of Hermite coefficients, see \cite{Vemuri2008hermite}. We prove that 
%% a certain combination of Hermite coefficients decays better than the 
%% one obtained in \cite{Vemuri2008hermite} for the individual Hermite 
%% coefficients.
%% \end{abstract}

%\begin{abstract}
%Assume that a function and its Fourier transform are dominated by 
%centered Gaussians of certain variance. Vemuri proved the decay of
%Hermite coefficients in terms of variance. In this article, we
%prove that a certain combination of Hermite coefficients decays better.
%\end{abstract}

\begin{abstract}
Assuming that a function and its Fourier transform are dominated by a
Gaussian, Vemuri found a sharp estimate for the decay rate of the
Hermite coefficients in terms of the variance of the dominating Gaussian.
Here we show that under the same assumptions, certain combinations
of Hermite coefficients have a better rate of decay.
\end{abstract}

\maketitle

%\tableofcontents\newpage

\section{Introduction}\label{S:intro}

If $f \in L^1(\R)$, the {\em Fourier transform} of $f$ is defined by
\begin{equation*}
\widehat f(\xi) = \frac{1}{\sqrt{2\pi}} \int f(x) e^{-i\xi x} \, dx.
\end{equation*}
Let $g_a(x) = e^{-ax^2/2}$.  In \cite{Garg2009}, the {\em Hardy Class} $H(a)$ is
defined by
\begin{equation*}
H(a) = \{ f \in L^1(\R) \st
f(x) = O\(g_a(x)\) \text{ and } \widehat f(\xi) = O\(g_a(\xi)\)\},
\qquad 0 < a < 1.
\end{equation*}
 
Let $\phi_n$ denote the $n$-th normalized Hermite function
(see \cite[Section 2.1]{ST}). Then $\{\phi_n\}_{n=0}^\infty$ forms
an orthonormal basis for $L^2(\R)$. It follows from the Cauchy-Schwarz
inequality and Mehler's formula (see \cite[Prop 2.3]{ST}) that if for
some $t>0$
\begin{equation}\label{E:exponential-decay}
\< f, \phi_n \> = O\(e^{-2nt}\), 
\end{equation}
then $f\in H(\tanh 2rt)$
%$f$ and $\widehat{f}$ are dominated by the Gaussian $e^{-\frac{1}{2}\tanh(2rt) x^2}$
%% \begin{equation}\label{E:Gaussian-domination}
%% \begin{aligned}
%%   f(x) &\; = O\(e^{-\frac{1}{2}\tanh(2rt) x^2}\), \qquad\text{and}\\
%%   \widehat{f}(x) &\; = O\(e^{-\frac{1}{2}\tanh(2rt) x^2}\)
%% \end{aligned}
%% \end{equation}
for $0<r<1$ (see e.g. \cite{Vemuri2008hermite}). Using other methods,
Radha and Thangavelu \cite{radha2009holomorphic} proved that under the
hypothesis (\ref{E:exponential-decay}), $f$ extends to $\mathbb{C}$ as an
entire function and satisfies the estimate
$$f(x+iy) = O\( e^{-\frac{1}{2} \tanh(2rt) x^2 +
  \frac{1}{2} \coth(2rt) y^2}\),$$
for $0<r<1$, and a similar estimate holds for
$\widehat{f}$ as well. 
%(see \cite{radha2009holomorphic}). 

%Suppose we only know that a function and its Fourier transform are 
%dominated by Gaussians. Then we would like to know if these conditions
%in turn imply some exponential decay of the Hermite coefficients of $f$. 
%A paper addressing similar question for the space of type W is published 
%by Janssen and Eijndhoven \cite{janssen1990spaces}.
On the other hand, Vemuri \cite{Vemuri2008hermite} proved that if $t>0$
and $f\in H(\tanh 2t)$ then
%% \begin{equation*}
%% \begin{aligned}
%%   f(x) &\; = O \(e^{-\frac{1}{2}\tanh(2t) x^2}\), \qquad\text{and}\\
%%   \widehat{f}(x) &\; = O\(e^{-\frac{1}{2}\tanh(2t) x^2}\)
%% \end{aligned}
%% \end{equation*}
%% then
\begin{equation}\label{E:Vemuri-estimate}
\< f, \phi_n \> = O\(n^{-1/4} e^{-nt}\),
\end{equation}
and that this estimate is sharp.

Shortly thereafter, Garg and Thangavelu generalized Vemuri's result
to the several variable case \cite{Garg2009}.
Combining the results of Vemuri with those of Radha and Thangavelu leads
to a certain loss of information about the decay properties of $f$ and
$\widehat{f}$.
Indeed, if we start with $f\in H(\tanh 2t)$, then by Vemuri's result
$\< f,\phi_n \> = O(n^{-1/4}e^{-nt})$, and so by Radha and Thangavelu's result,
we get $f\in H(\tanh rt)$ for $0<r<1$.
Therefore, it stands to reason that the condition $f\in H(\tanh2t)$
has consequences on the Hermite coefficients of $f$ other than mere
exponential decay. Indeed, our main result is a sort of pair correlation
between certain Hermite coefficients.
%%In this paper, we obtain an improved bound for certain
%%combinations of Hermite coefficients. Our main result is the following.

\begin{thm}\label{T:resonance}
Let $t>0$. If $f\in H(\tanh 2t)$ then
\begin{equation*}
\<f,\phi_n\>+ \frac{n(n+2)}{\sqrt{\prod\limits_{j=1}
^{4}(n+j)}} e^{4t} \<f,\phi_{n+4}\>= O\(n^{-3/4} e^{-nt}\),
\end{equation*}
for $n=1,2,\dots$.
\end{thm}
\begin{rem}
In \cite{Vemuri2008hermite}, Vemuri used (\ref{E:Vemuri-estimate}) to prove
a uniform Gaussian bound for the solution of the harmonic oscillator Schr\"{o}dinger 
equation with initial solution $\phi_0\in H(\tanh 2t)$, and conjectured the sharp bound.
Recently, substantial progress towards the conjecture was 
made in \cite{kulikov2023gaussian}.  We hope that the Theorem \ref{T:resonance} could
be useful in this direction.
\end{rem}

%% \begin{rem}\label{R:betterdecay}
%% There are functions in $H(a)$ for $a\in(0,1)$, whose above mentioned
%% combination of Hermite coefficients decays better  than the one
%% obtained in the Theorem \ref{T:resonance}. We provide an example of
%% such function at the end of this paper.
%% \end{rem}

%% \begin{rem}\label{R:comparison}
%%  If $f\in H(a)$ then the bound obtained in Theorem \ref{T:resonance}
%%  for the mentioned combinations have a better decay than the bound for
%%  individual Hermite coefficients. 
%% \end{rem}

%% In \cite{Vemuri2008hermite}, the main idea to prove (\ref{E:Vemuri-estimate})
%% is to obtain estimates of the Bargmann transform (see \cite [Section 2.1]{ST})
%% and improve these estimates by applying Phragm\'en-Lindel\"of principle. 
%% Here, to prove Theorem \ref{T:resonance}, we control the bad growth of the
%% improved estimate of the Bargmann transform by the use of certain polynomial.

%% There are two steps in the proof of (\ref{E:Vemuri-estimate}) in
%% \cite{Vemuri2008hermite}. The first step is to use the hypotheses
%% on $f$ together with the Phragmen-Lindelof principle to estimate
%% the growth of the Bargmann transform $Uf$ of $f$.
%% The second step is to use this growth estimate together with
%% Cauchy's integral formula along a well chosen contour to estmate
%% the Taylor coefficents of $Uf$ at $0$; the latter are closely
%% related to the Hermite coefficients of $f$ (see \cite{Vemuri2008hermite}).
%% In this work, we use we modify the second step to obtain better estimates for
%% certain combinations of Taylor coefficients of $Uf$. 

\section{The Proof}

We first recall some properties of the Bargmann transform.
Let $\H$ denote the Hilbert space of all entire functions $F$ on $\C$
such that
\begin{equation*}
\norm{F}^2_\H = \int_{\C} \abs{F(w)}^2 \,
                        \frac{e^{-\abs{w}^2/2} \, du \, dv}{\sqrt{4\pi}}
             < \infty \quad\text{($w=u+iv$)}.
\end{equation*}
For a Schwartz class function $f$, the {\em Bargmann transform} of $f$
is defined by
\begin{equation*}
Bf(w) = \frac{e^{-w^2/4}}{\sqrt{\pi}} \int_{\R} e^{xw} e^{-x^2/2} f(x) \,dx.
\end{equation*}
It is shown in \cite[Section 2.1]{ST} that $B$ extends to an isometric
isomophism
from $L^2(\R)$ to $\H$, and
\begin{equation*}
B\phi_n(w) = \frac{w^n}{\sqrt{2^n n!\pi^{1/2}}}.
\end{equation*}
Let $a\in (0,1)$. Let $w=re^{i\theta}$ and $\mu=\frac{1-a}{1+a}$. It is
shown in \cite{Vemuri2008hermite} that if $f \in H(a)$ then
\begin{equation}\label{E:B1}
\abs{Bf(w)}\le C \sqrt{\frac{2}{1+a}}
	\exp \frac{(\mu + (1-\mu)\sin^2 \theta)r^2}{4},
\end{equation}
\begin{equation}\label{E:B2}
\abs{Bf(w)}\le C \sqrt{\frac{2}{1+a}}
	\exp \frac{(\mu + (1-\mu)\cos^2 \theta)r^2}{4}.
\end{equation}
Furthermore, it is shown, using the Phragm\'en-Lindel\"of principle
that
\begin{equation}\label{E:B3}
\abs{Bf(w)} \le C \sqrt{\frac{2}{1+a}}
\exp\(\sqrt{\mu} \sin\(\abs{2\theta}\)\frac{r^2}{4}\), 
\end{equation}
for $\theta_0 \le \theta-\frac{k\pi}{2}\le \frac{\pi}{2}-\theta_0$,
$k=0,1,2,3$, where $\theta_0=\tan^{-1}(\sqrt\mu)$.
Observe that $0<\theta_0<\frac{\pi}{4}$.
 Now, define $\gamma_n(t)=\(\frac{4n(n+2)}{\mu}\)^{1/4}e^{it}$
 for $0\le t\le2\pi$. 
By the Cauchy integral formula for derivatives, we have $Bf(w)=\sum_
{n=0}^{\infty}c_n w^n$ where
\begin{equation*}
c_n
= \frac{1}{2\pi i} \int_{\gamma_n} \frac{Bf(w)}{w^{n+1}} \,dw.
\end{equation*}
Therefore
\begin{equation*}
\abs{c_n+\frac{4n(n+2)}{\mu}c_{n+4}}=\frac{1}{2\pi}\abs{\int_{\gamma_n}
\(w^4+\frac{4n(n+2)}{\mu}\)\frac{Bf(w)}{w^{n+5}}dw}.
\end{equation*}
Thus
\begin{equation}\label{E:dihedral}
\begin{aligned}
\abs{c_n+\frac{4n(n+2)}{\mu}c_{n+4}}
                   \le&\; \frac{1}{2\pi}\(\frac{\mu}{4n(n+2)}\)^{n/4}
                          \int_{0}^{2\pi}
                          \abs{e^{4it}+1}\abs{Bf\(\(\frac{4n(n+2)}{\mu}\)^{1/4}
                           e^{it}\)}dt\\
                    =&\; \frac{1}{\pi}\(\frac{\mu}{4n(n+2)}\)^{n/4} \sum_{k=1}^4 
                         \int_{\frac{(k-1)\pi}{2}}^{\frac{k\pi}{2}}
                         \abs{\cos2t}\abs{Bf\(\(\frac{4n(n+2)}{\mu}\)^{1/4}
                           e^{it}\)} \,dt.
\end{aligned}
\end{equation}
%% By Lemma \ref{L:criticalpts}, we have
%% \begin{equation}\label{E:dihedral}
%% \abs{c_n+\frac{4n(n+2)}{\mu}c_{n+4}}\le 
%% \frac{1}{\pi}\(\frac{\mu}{4n(n+2)}\)^{n/4} 
%% \sum_{k=1}^4 \int_{\frac{(k-1)\pi}{2}}^{\frac{k\pi}{2}}
%% \abs{\cos2t}\abs{Bf\(\(\frac{4n(n+2)}{\mu}\)^{1/4}e^{it}\)} \,dt.
%% \end{equation}
By inequalities (\ref{E:B1}), (\ref{E:B2}),  and (\ref{E:B3}) we have
\begin{equation*}
\int_{0}^{\pi/2}\abs{\cos2t}\abs{Bf\( \(\frac{4n(n+2)}{\mu}\)^{1/4} e^{it}\)}\, dt
\le C\sqrt{\frac{2}{1+a}} (I_n+J_n+K_n)
\end{equation*}
where
\begin{equation*}
\begin{aligned}
I_n
=&\; \int_0^{\theta_0} \abs{\cos2t} \exp\(\frac{\sqrt{n(n+2)}
     (\mu+(1-\mu)\sin^2t)\, }
        {2\sqrt{\mu}}\) dt,\\
J_n
=&\; \int_{\theta_0}^{\frac{\pi}{2}-\theta_0}\abs{\cos2t}
     \exp\(\frac{\sqrt{n(n+2)}}{2}\sin2t\) \, dt,
         \quad\text{and} \\
K_n
=&\; \int_{\frac{\pi}{2}-\theta_0}^{\frac{\pi}{2}}\abs{\cos2t}
     \exp\(\frac{\sqrt{n(n+2)}(\mu+(1-\mu)
         \cos^2t)}{2\sqrt{\mu}}\)\, dt.
\end{aligned}
\end{equation*}
For $n\in\N$ let $h_n:(\theta_0,\frac{\pi}{4})\rightarrow\R$ be defined by
\begin{equation*}
h_n(t)=\sin2t+\frac{2}{\sqrt{n(n+2)}}\log \abs{\cos2t}.
\end{equation*}
Clearly, $h_n(t) = h_n(\frac{\pi}{2} - t)$ {\em for all} $t\in(\frac{\pi}{4},\frac{\pi}{2}-\theta_0)$.
It follows that
\begin{equation}\label{E:Jn}
J_n
= \int_{\theta_0}^{\frac{\pi}{2}-\theta_0} 
     \exp\(\frac{\sqrt{n(n+2)}}{2} h_n(t)\) \, dt 
 = 2\int_{\theta_0}^{\frac{\pi}{4}}\exp\(\frac
          {\sqrt{n(n+2)}}{2} h_n(t) \)dt.
\end{equation}

Since
\begin{equation*}
\(\sqrt{\frac{n}{n+2}}\)^2 + \(\sqrt{\frac{2}{n+2}}\)^2 = 1,
\end{equation*}
there exists a unique $t_0\in(0,\frac{\pi}{4})$ such that
\begin{equation*}
\sin2t_0 = \sqrt{\frac{n}{n+2}} \quad \text{and} \quad \cos2t_0 = \sqrt{\frac{2}{n+2}}.
\end{equation*}
Thus, we have
\begin{equation}\label{E:hn'}
h_n'(t_0) = 2\cos2t_0 - \frac{4\tan2t_0}{\sqrt{n(n+2)}} = 0,
\end{equation}
and
\begin{equation}\label{E:hn''}
h_n''(t_0) = -4\sin2t_0 - \frac{4\sec^2 2t_0}{\sqrt{n(n+2)}} = 
\frac{-8(n+1)}{\sqrt{n(n+2)}} < 0.
\end{equation}
Also, observe that, for large enough 
$n\in\N$, we have $\theta_0 \le t_0 < \frac{\pi}{4}$.
Therefore, we may estimate
$J_n$ by the use of Laplace's method (see \cite[Section 2.4]{Erdelyi}), 
which can be stated as follows.
\begin{thm}\label{T:Laplace}
Let $G, H : [\alpha, \beta] \rightarrow \R$ be continuous and twice continuously 
differentiable functions. If $H'(t_0) = 0$, and $H''(t_0) < 0$, for a
unique $t_0 \in [\alpha, \beta]$. Then  
\begin{equation*}
\int_{\alpha}^{\beta} G(t) e^{xH(t)} dt \sim G(t_0) e^{xH(t_0)} 
\[\frac{-\pi}{2xH''(t_0)}\]^{1/2}.
\end{equation*}
\end{thm}
In our case (see equation (\ref{E:Jn})), $G(t) = 1$, $x = \frac{\sqrt{n(n+2)}}{2}$, 
and $H(t) = h_n(t)$.
Therefore
\begin{equation}\label{E:Laplace1}
\[\frac{-\pi}{2xH''(t_0)}\]^{1/2} = \[\frac{\pi}{8(n+1)}\]^{1/2},
\end{equation}
and
\begin{equation}\label{E:Laplace2}
e^{xH(t_0)} = e^{\frac{\sqrt{n(n+2)}}{2}\(\sqrt{\frac{n}{n+2}} + 
\frac{2}{\sqrt{n(n+2)}} \log \sqrt{\frac{2}{n+2}}\)} = \sqrt{\frac{2}{n+2}}e^{n/2}.
\end{equation}
From equations (\ref{E:Jn}), (\ref{E:hn'}), (\ref{E:hn''}), (\ref{E:Laplace1}), 
(\ref{E:Laplace2}), and Theorem \ref{T:Laplace}, we get
\begin{equation*}
J_n \sim \sqrt{\pi}\, n^{-1}e^{n/2}.
\end{equation*}

Observe that
\begin{equation*}
\mu+(1-\mu)\sin^2t\le\mu+(1-\mu)\sin^2\theta_0 \quad \text{\em for all}
\quad t\in\[0,\theta_0\],
\end{equation*}
and
\begin{equation*}
\mu+(1-\mu)\sin^2\theta_0\le \sqrt{\mu}\sin2t\quad \text{\em for all} 
\quad t\in\[\theta_0,\frac{\pi}{4}\].
\end{equation*}
Therefore
\begin{equation*}
I_n \le \sin2\theta_0 \exp\(\frac{\sqrt{n(n+2)}(\mu+(1-\mu)\sin^2\theta_0)\, }
        {2\sqrt{\mu}}\),
\end{equation*}
and
\begin{equation*}
J_n \ge 2(1-\sin2\theta_0) \exp\(\frac{\sqrt{n(n+2)} (\mu+(1-\mu)\sin^2\theta_0)\, }
{2\sqrt{\mu}}\).
\end{equation*}
%Therefore
%\begin{equation*}
%I_n\le \(\frac{2\sin2\theta_0}{1-\sin2\theta_0}\) \int_{\theta_0}^{\frac{\pi}{4}}
%\cos2t \exp\(\frac{\sqrt{n(n+2)}}{2}\sin2t\) \, dt.
%\end{equation*}
Thus
\begin{equation*}
  I_n = O\(J_n\).
\end{equation*}
Clearly
\begin{equation*}
  K_n = I_n.
  \end{equation*}
%% Similarly, we observe that
%% \begin{equation*}
%% \mu+(1-\mu)\cos^2t\le\mu+(1-\mu)\cos^2\theta_1 \quad \text{\em for all}
%% \quad t\in\[\theta_1,\frac{\pi}{2}\],
%% \end{equation*}
%% and
%% \begin{equation*}
%% \mu+(1-\mu)\cos^2\theta_1\le \sqrt{\mu}\sin2t\quad \text{\em for all} 
%% \quad t\in\[\frac{\pi}{4},\theta_1\].
%% \end{equation*}
%% Therefore
%% \begin{equation*}
%% K_n\le \(\frac{\sin2\theta_1}{1-\sin2\theta_1}\) \int_{\frac{\pi}{4}}^{\theta_1}
%% \cos2t \exp\(\frac{\sqrt{n(n+2)}}{2}\sin2t\) \, dt.
%% \end{equation*}
%% Since $\theta_1=\frac{\pi}{2}-\theta_0$, therefore we have
%% \begin{equation*}
%% I_n+K_n\le \(\frac{\sin2\theta_0}{1-\sin2\theta_0}\)J_n.
%% \end{equation*}
Hence, we get
\begin{equation*}
\int_{0}^{\pi/2}\abs{\cos2t}\abs{Bf\( \(\frac{4n(n+2)}{\mu}\)^{1/4} e^{it}\)}\, dt
=O\(n^{-1}e^{n/2}\).
\end{equation*}
The other three integrals in (\ref{E:dihedral}) are also
$O\(n^{-1}e^{n/2}\)$ by equation (\ref{E:B3}),
and the fact that the right hand sides of inequalities (\ref{E:B1})
and (\ref{E:B2}) do not change when we replace
$\theta$ by $\pi-\theta$ or $2\pi-\theta$.
We conclude from equation (\ref{E:dihedral}) that
\begin{equation}\label{E:Taylor's-estimate}
c_n+\frac{4n(n+2)}{\mu}c_{n+4}=O\[n^{-1}\(\frac{\sqrt{\mu}e}{2n}\)^{n/2}\].
\end{equation}

Since
\begin{equation*}
\begin{aligned}
\<f, \phi_n\>
  =&\; \<Bf, B\phi_n\>\\
  =&\; \int\int \(\sum_{k=0}^\infty c_k w^k\)
                \bar{\(\frac{w^n}{\sqrt{2^n n!\pi^{1/2}}}\)}
                \frac{e^{-r^2/2} \, du \, dv}{\sqrt{4\pi}}\\
  =&\; \frac{c_n}{\sqrt{2^n n!\pi^{1/2}}}
       \int\int r^{2n} \frac{e^{-r^2/2} \, du \, dv}{\sqrt{4\pi}}\\
  =&\; \sqrt{2^n n!\pi^{1/2}} c_n.
\end{aligned}
\end{equation*}
Therefore, we obtain
\begin{equation*}
\<f,\phi_n\>+\frac{1}{\mu}\frac{n(n+2)}{\sqrt{\prod\limits_{j=1}^{4}
(n+j)}}\<f,\phi_{n+4}\>= \sqrt{2^n n!\pi^{1/2}}\[c_n+\frac{4n(n+2)}{\mu}
c_{n+4}\].
\end{equation*}
It follows from equation (\ref{E:Taylor's-estimate}) and Stirling's
formula that
\begin{equation*}
\<f,\phi_n\>+\frac{1}{\mu}\frac{n(n+2)}{\sqrt{\prod\limits_{j=1}^{4}
(n+j)}}\<f,\phi_{n+4}\>=O\(n^{-3/4}\ \mu^{n/4}\).
\end{equation*}

Taking $\mu = e^{-4t}$ gives the result stated in Theorem \ref{T:resonance}.

%% Now we give an example of a function in $H(a)$ for $a\in(0,1)$,
%% whose above mentioned combination of  Hermite coefficients have
%% better decay. Let
%% \begin{equation*}
%% f(x)=\exp\(\frac{-a+i\sqrt{1-a^2}}{2}x^2\)
%% \end{equation*}
%% where $a\in(0,1)$. Observe that $f\in H(a)$. However,
%% \begin{equation*}
%% \begin{aligned}
%% Bf(w)
%% =&\; \frac{e^{-w^2/4}}{\sqrt{\pi}}  \int e^{xw} e^{-x^2/2}
%%      \exp\(\frac{-a+i\sqrt{1-a^2}}{2}x^2\) \, dx\\
%% = & \; \frac{1}{\sqrt{\pi}} \exp\(-i\sqrt{\mu}\frac{w^2}{4}\)\\
%% = & \; \frac{1}{\sqrt{\pi}} \sum_{n=0}^{\infty} \frac{(-i)^n \mu^{n/2}
%%           w^{2n}}{4^n n!}.
%% \end{aligned}
%% \end{equation*}
%% Therefore $\langle f, \phi_n\rangle=0$ for $n$ odd and
%% \begin{equation*}
%%   \<f, \phi_n\> = \frac{(-i)^{n/2}2^{-n/2}\mu^{n/4}\sqrt{n!}}{\sqrt{\pi}
%%                   (\frac{n}{2})!},
%% \quad n=0,2,4,\dots.
%% \end{equation*}
%% Thus, by Stirling's formula, we have
%% \begin{equation*}
%% \<f,\phi_n\>+\frac{1}{\mu}\frac{n(n+2)}{\sqrt{\prod\limits_{j=1}^{4}(n+j)}}
%% \<f,\phi_{n+4}\>\sim\(\frac{2}{\pi^3}\)^{1/4}\[n^{-5/4}\(\frac{1-a}{1+a}\)^{n/4}\] 
%% \quad n=0,2,4,\dots.
%% \end{equation*}

\section*{Acknowledgements}
I would like to thank my advisor M. K. Vemuri for suggesting this
problem to me.

\section*{Declarations}
I hereby declare that this work has no conflict of interest-neither
personal nor financial.

\bibliographystyle{amsplain}
\bibliography{v7-resonance}

\end{document}